# Alternative to the Well-known Statistical Dynamics of Linear Systems

V. N. Tibabishev

Private laboratory identification of dynamic systems, office 151/36, Lugansk, 91004, Ukraine

**Abstract**   The problem of determining the mathematical model of the dynamics of multi-dimensional control systems in the presence of noise under the condition that the correlation functions cannot be found. Known statistical dynamics of linear systems is a more effective alternative. Background information is presented in the form of individual implementations nonergodic stochastic processes. Such a realization is deterministic functions. We introduce the concept of systems of sets of signals for the components on the semiring. For the system of sets of linearly dependent and linearly independent of the measured signals of a certain frequency properties. Frequency method is designed to deal with the noise on the set of deterministic functions. Example is the determination of the dynamic characteristics of the aircraft in accordance with the data obtained in one automatic landing.

**Keywords**   Semiring deterministic components of the signals, The Hilbert space of almost periodic functions, Anti-interference

## 1. Introduction

Many applied problems of control systems were solved by the statistical dynamics in the presence of uncorrelated noise at present. In the theory of stationary random process occurs dualism. The same random stationary process has two incompatible models of representation. One submission random processes have a continuous spectrum, belong to the Hilbert space $L_2(-\infty, +\infty)$ by Wiener and have the ergodic property[1]. According to another model, random stationary process has a discrete spectrum, belongs to the Hilbert space of almost periodic functions $B_2(-\infty, +\infty)$ by E. Slutsky[2] and does not have the ergodic property[1]. A widespread model for Wiener thanks ergodic property. Primary data are realizations of random processes in the presence of interference. Uncorrelated noise is suppressed by correlation analysis. Correlation functions are secondary source data.

It is known[3] that any vector can have up to a countable set of non-zero projections for any orthonormal system of power in any Hilbert space. It follows that the Hilbert spaces there is no process of continuous spectrum. The notion of a continuous spectrum in the Hilbert space $L_2(-\infty, +\infty)$ entered incorrectly due to an incorrect assessment of the domain of the Fourier - Plancherel[4].

Stationary random processes can have only discrete spectrum. Stationary random processes with discrete ergodic have first order expectation and have no ergodicity second order dispersion and correlation functions[1]. Secondary raw data can only be found by averaging over the set of recorded simultaneously under identical environmental conditions. The implementation of this method of obtaining secondary source of data is limited. First, in many cases, impossible to get a lot of recorded simultaneously.[5] Second, if it's there, it is usually impossible to obtain multiple implementations for the same external conditions. For example, it is theoretically possible to provide multiple synchronized landing aircraft on many runways (runway). Put a lot of runway in one location on earth is almost impossible. Extreme runway can be located far away from each other. Recorded implementation cannot be obtained simultaneously in the same weather conditions.

Identification of dynamic characteristics we consider a specific example. The mathematical model of the dynamics of the aircraft is difficult to obtain analytically because of elastic deformation airframe. Therefore, the dynamic characteristics found by solving the problem of identification of the experimental data obtained in single automatic landing an airplane class of IL-96 on a runway (runway). The plane is a multidimensional object of control. The dynamics of the aircraft should be determined only on the two control channels. The first control channel is given pitch - x_1 - actual pitch - y in the electric remote control system. The second channel is the position of the control levers engine (ORE) - x_2 - pitch plane - y. On the pitch the aircraft affect the angular position of the flaps x_3 and position angle of the stabilizer x_4. The angular position of the flaps and the angle of the stabilizer change the motion of the aircraft on final approach. Therefore the dynamics of the aircraft and the output are no stationary. In single auto-land were obtained synchronous recording input $x_i$, i = 1,2,3,4, and the output y from 274 reference in 0.5 seconds one copy.

The well-known method for solving the problem in the class of linear time-varying systems is reduced to solving a system of integral equations of the first kind[6].

$$\sum_{i=1}^{d} \int_{t-T}^{t} R_{x1xi}(\theta,\tau) k_{1i}(t,\tau) d\tau = R_{y1x1}(t,\theta); \quad (1)$$
$$\sum_{i=1}^{d} \int_{t-T}^{t} R_{xdxi}(\theta,\tau) k_{di}(t,\tau) d\tau = R_{ydxd}(t,\theta),$$

where d - the number of inputs and outputs.

The system of integral equations contains all the autocorrelation function of the input signal $R_{xixi}$ all the cross-correlation function between the input signals $R_{xpxi}$ $\forall p \neq \forall i$ and all the cross-correlation function between all input and output signals $R_{ypxi}$. It is known[7] that the non-stationary random processes do not have the ergodic property. Therefore, the correlation function is obtained by averaging the correct time unit implementations possible. Contracting authority may provide the raw data only once.

In addition to get the system of equations (1) to four inputs have four independent outputs. In this case there is only one way out. This problem cannot be solved by a known method for these reasons.

It is known[7] that if the implementation of a random process is defined as a single implementation, it is a normal (non-random) function. Thus, the problem of solving the problem of identification of dynamic characteristics of a control in the presence of noise is not on a set of random functions, and the simultaneous recording of input and output signals determined to be submitted in one copy.

The plane is nonlinear and non-stationary object of control.[8] However, in a small interval of observation with an almost constant speed of flight on the final approach the plane can be considered linear and stationary control. Therefore, we will solve the problem of identifying the dynamic characteristics of the first approximation in the class of linear time-invariant models.

The task of identifying the dynamic characteristics of the solution precedes suppression. This problem can be solved if the useful signal and noise differ in their properties. Find it necessary to construct mathematical models of signal and noise in the class of non-random functions on a single specimen.

## 2. Processes and Signals are Measurement

All that is described as a function of time is called a process. Process, converted into electrical form sensor (analog or digital) suitable for further processing, we call signal. The exact component of the signal is distorted by additive and multiplicative noise. Management systems may have one or more control channels. Control channel - a device displaying many precise components of the input signals to (in) the set of precise components of the output signals. By accurate data components of the measured signals are understood as functions that are associated operator equation. Mathematical model of the control channel is defined by the operator. Operators can be linear and nonlinear stationary and non-stationary. Multivariate control system is divided into systems with linearly independent and linearly dependent inputs. We consider multidimensional systems with linearly independent inputs initially.

In the simplest case of one-dimensional dynamics of the control object in a first approximation, is in a class of linear stationary models, which describe the differentiation operator with constant coefficients, for example n- the first order $Dy = x$. Inverse operator of differentiation is an integral operator of convolution type $Ax = y$. The output signal y is the image of the input signal x, which is called a prototype. It will be shown that for linear operators and inverse images are linearly dependent processes.

Precise input process x by the sensor is converted into an electrical signal, which is measured by the input signal. As a result of the conversion and transmission of communication channels accurate signal distortion additive and multiplicative noise $\tilde{x} = \vartheta x + n$, where $\vartheta$ and n - respectively multiplicative and additive noise. These disturbances are of different nature education and are therefore linearly independent processes.

Instead of just the output y is available only measured output $\tilde{y} = \theta y + m$, $\theta$ and m - respectively multiplicative and additive noise output. Additive and multiplicative noise generated at the output not only sensors, but also disturbances. For example, in the form of wind gusts upward and downward flow of air masses that affects the output signals of the aircraft. Nature Education noise $\vartheta$ and n is different from the nature of the formation of interference $\theta$ and m. On this basis, bude assume that the processes $\vartheta$, n, $\theta$ and m are linearly independent processes.

It is known[9] that the multiplicative interference can lead to additive noise generated by multiplicative noise. The observed input to the form

$$\tilde{x} = \hat{\vartheta} x + p + n, \quad (2)$$

where $\hat{\vartheta}$ - a number equal to the mean value of the multiplicative noise $\vartheta$,

p -reduced additive noise generated by the multiplicative interference $\vartheta$ .

Additive noise $p$ and n do not generate the forced motion control channel. They appear in the measurement channel. Precise input $x$ process and additive noise $n$ and $p$ have different nature of education, so are pairwise linearly independent processes.

The observed output to the form $\tilde{y} = \hat{\theta}y + \gamma + m$, where $\hat{\theta}$ - a number equal to the mean value of the multiplicative noise $\theta$, $\gamma$ — reduced additive noise generated by the multiplicative interference $\theta$.

Additive noise $\gamma$ and $m$ are not part of the forced movement. They have a different nature arise. Therefore, the exact process $y$ and additive noise $\gamma$ and m are pairwise linearly independent processes.

We shall consider the control system, the dynamics of which is described by the Newton[10], according to which the mass of the body $m$, the second derivative of the function $\ddot{x}$, for example, $x$ describes the progressive movement is related to the power $\gamma$ of ordinary inhomogeneous differential equation $m\ddot{x} = \gamma$. The aircraft move on final approach to weight almost constant over time. Newton differential equation is linear stationary equation. Therefore, under the symbol $x$ can be any component of the signal. For each component of the measured signal (2), for example, n condition for the existence and uniqueness of solutions of differential equations of motion $m\ddot{n} = \gamma_n$ m performed. It follows that the second derivative $\gamma_n$ and $\ddot{n}$ is a continuous bounded function[11]. Obviously, the first derivative $\dot{n}$ of the function itself $n$ is also a continuous function.

It is known[11] that in this case the differential equation will be presented convergent trigonometric series. If the function $\gamma$ is a continuous function on the whole line, it just seems to be convergent trigonometric series[12]. Furthermore, we assume that the elements of the measured processes belong to the set of bounded continuous functions representing the convergent trigonometric series. On the set of functions square, for example, of the input is not integrable on the entire line of Riemann and Lebesgue-Stieltjes integrable and the function $g(t) = t/2T$
$$F = \lim_{T\to\infty} \int_{-T}^{+T} x^2(t) dg(t) = \lim T \to \infty 1 2 T - T + T x 2 t d t \equiv M x 2 (t) < \infty. \quad (3)$$

These properties are, for example, the functions that belong to the Hilbert space of almost periodic functions $B_2(-\infty, +\infty)$ in the sense of Besicovitch.[13] Furthermore, we assume that the components of the observed input and output signals are elements of the Hilbert space of almost periodic functions $B_2(-\infty, +\infty)$ in the sense of Besicovitch.

## 3. System of Linearly Independent Sets of Signals

We now consider the multidimensional linear systems. Let control system has d inputs and h outputs. System sets are highlighted by several sets of measured signals. For example, the input signal is measured k $\tilde{x}_k = \hat{\theta}x_k + p_k + n_k$ belongs to a system of sets $\tilde{x}_k \in \mathbb{R}_k = X_k \cup P_k \cup N_k$, where $\hat{\theta}x_k \in X_k$ -many precise components of the signals that generate forced motions in the control channels. $P_k$ is reduced additive noise generated by the multiplicative interference $\theta$, where $\hat{\theta}$ - the average. $N_k$-set of additive noise, k -dead entrance. Another system of sets of linearly independent elements is the union of all input signals $\mathbb{R} = \cup_{k=1}^{k=d} \mathbb{R}_k$. Similarly, we obtain a system of linearly independent sets of output signals $\mathbb{Z} = \cup_{k=1}^{k=h} \mathbb{Z}_k$, where $\tilde{y}_k \in \mathbb{Z}_k = Y_k \cup \Gamma_k \cup M_k$, and $\tilde{y}_k = \hat{\eta}y_k + \gamma_k + m_k$. $Y_k$ - many precise output. $\Gamma_k$ - set of reduced additive noise, generated by the multiplicative interference $\gamma_k$, where $\hat{\eta}$ - the average. $M_k$ - set of additive noise, k -died out. This system of linearly independent sets of elements observed at the output, you can expand it by adding a system of sets of linearly independent components of the noise, distorting the input signals.

The set of elements, such as, $X_k$ conceived as a whole, in the sense that each element of the set is the linear span of a certain orthonormal basis of harmonic $f_{kl} \in G_k$, $l = 1,2,3,...,f_{kl} = \exp(j\omega_{kl}t)$. Moreover, each i - th actual element, for example, the exact component of input signal at the input k is represented as
$$x_{ki} = \sum_{l=1}^{l=c}(a_{kil}\cos(\omega_{kl}t) + b_{kil}\sin(\omega_{kl}t)). \quad (4)$$

If the coefficients $a_{kil}$ and $b_{kil}$ are uncorrelated with each other random numbers, depending on the index implementation i, then the elements $x_{ki}$ are wide-sense stationary stochastic processes that, when averaged over the set of realizations generate the correlation function of the general form $R_k(\tau) = \sum_l A_{kl} \cos(\omega_{kl}t)$ [2].

Along with the presentation of a stationary random process in the form of (4) there are continuous stationary random processes with orthogonal increments[1]. Such random processes do not have a derivative at any point. Apparently, not differentiable random processes describe the processes, the dynamics of which is not subject to classical Newtonian mechanics, and beyond that we have adopted restrictions.

The functional spaces can the sets notion of norm for the element. Observed signals are the union of different properties of sets of precise components and noise. In general, each component of the signal is closely accurate observed in the presence of two additive noise belonging to different sets. For the elements of the input set the measure $\mu$ (V) - numerical characteristic that satisfies the conditions
1. the domain of the function $\mu$ (V) is semi-ring sets
2. value $\mu$ (V) of are real and non-negative,
3. $\mu$ (V) is additive, for any finite decomposition $V = V_1 \cup ... \cup V_k$ into disjoint sets the equality $\mu(V) = \sum_{l=1}^{c} \mu(V_l)$[14].

It is known[14] that, as a measure can be taken, for example, the integral of a non-negative function. We take as a measure of the value sets of the measured signals of the functional (3), the domain of the semi-circle which must be set. Set-system functions $\mathbb{R}$ is a semiring if it is the union of $\mathbb{R} = \cup_{k=1}^{k=d} R_k$ disjoint sets satisfying $R_k \cap R_m = \emptyset$ for all $k \neq m$, where $\emptyset$- empty 'symbol sets.

Let $x_{ki} \in R_k$ arbitrarily chosen element of the set $R_k$, $1 \leq k \leq d$, $i = 1,2,3...$, $\tilde{x}_{ki} = \hat{\vartheta}x_{xki} + p_{xki} + n_{xki}$. Functional (3), defined for randomly selected items from the union of sets $\mathbb{R} = \cup_{k=1}^{k=d} R_k$, $F(V) = M\{\sum_{k=1}^{k=d} \tilde{x}_{ki}\}^2$ will be additive numerical characteristic if all inner products to randomly selected sets of components $(\tilde{x}_{ki}, \bar{\tilde{x}}_{mi}) \equiv M\{\tilde{x}_{ki} \cdot \bar{\tilde{x}}_{mi}\} = 0$ for all $k \neq m$, where $\tilde{x}_{ki}$ and $\tilde{x}_{mi}$ - complex conjugate functions. It is known[14] that if the elements are linearly independent, then they are orthogonal. In determining the functional (3) all scalar products are pairwise linearly independent, and hence the orthogonal functions. It follows that the elements are linearly independent systems of sets in Hilbert space of almost periodic functions is the semicircle, and the functional (3) defines a measure on it.

Check the condition $R_k \cap R_m = \emptyset$ on infinite sets of processes is almost impossible. Let $V_k, V_m$ are metric spaces and contain n elements. For each element $x_i \in V_k$ and every item $y_l \in V_m$ can define $n^2$ values of the distances $\rho(x_i, y_l)$. If for all $i, l = 1,2,3,...n$ the condition $\rho(x_i, y_l) \neq 0$, then the suppression $V_k \cap V_m$ is empty. Such an algorithm cannot be implemented in practice, for the reason that each set $R_k$ contains an infinite countable set of elements $n = \infty$. In this connection there is a need to develop a method implemented in practice to determine the half-ring of infinite sets of components of the observed signals.

Theorem 1. Given a system of infinite sets of functions $\mathbb{R} = \cup_{k=1}^{k=q} V_k$, where q - the number of components of the measured signal is spanned $V_k = \{x_{ik} \in V_k: x_{ik} = \sum_p a_{pki} e_{pk}\}$ spanned by the corresponding orthonormal bases $e_{pk} = exp(j\omega_{pk}t)$, $p, i = 1,2,3,...$ in the Hilbert space of almost periodic functions. If the system is a set of functions $\mathbb{R}$ semicircle, the numerical set of angular frequencies $\omega_{pk} \in \Omega_k$ harmonic bases, forming a system of infinite sets of functions of $\mathbb{R}$, is a semicircle of disjoint sets of numbers of the angular frequency harmonics orthonormal bases $\Omega_k \cap \Omega_l = \emptyset$ for all k ≠ l.

Proof of Theorem 1. Since R is a semi-circle, the functional (3) will determine the additive measure, if all $x_{ki} \in V_k, x_{pl} \in V_p$, $p, i, l = 1,2,3,..., k = 1,2,3,..q$ q scalar products are zero for all $k \neq p$. The scalar product is zero for all $a_{ki} \bar{b}_{pl}$ at randomly selected i and l, if for all $k \neq p$ is the condition $\omega_k - \omega_p \neq 0$

$$M\{\sum_i^c a_{ki} exp(j\omega_k t) \sum_l^c \bar{b}_{pl} exp(-j\omega_p t)\} =$$
$$M\{\sum_i^c \sum_l^c a_{ki} \bar{b}_{pl} exp j(\omega_k - \omega_p)t\} = 0. \quad (5)$$

It follows that the system of numerical sets of angular frequency $\mathbb{Q} = \cup_{k=1}^{k=d} \mathbb{Q}_k$ a semicircle where $\mathbb{Q}_k = \cup_{kc=1}^{c=q} \Omega_{kc}$_kc, $k$ - number of input, $q$ - the number of independent components of the signal.

Consequence of Theorem 1. It is known[3] that in the Hilbert space of almost periodic functions defined by a completely continuous symmetric operator normal form of convolution type, displaying $B_2(-\infty, +\infty)$ to itself. It is known[15] that if A is a compact and normal form, it generates only a finite number of basis functions and the non-zero eigenvalues. Therefore, in terms of elements (4), and the scalar products (5) the value of $c < \infty$. It follows that the system $\mathbb{Q} = \cup_{k=1}^{k=d} \Omega_k$ number sets values of the angular frequency of the harmonic basis functions of independent components of the input signals is the union of pairwise disjoint finite sets of numbers.

Corollary 2 of Theorem 1. It is known[2] that the set of almost periodic functions of random stationary process is a trigonometric series (4) with real numbers $\omega_k$. In the theory of random processes, the domain of the series (4) is just the set of almost periodic functions. The components of the measured signal as well represented by the series (4), but their domain is not a separate set of components, and the system is set to be a semicircle. In this case, the frequency $\omega_k$ choice additional restrictions. Batch processes are represented by harmonic series (4) with multiple frequencies also belong to the set of almost periodic functions. For example, suppose the set of frequencies $\Omega_a$ contains frequencies $\omega_k = k\omega_1$, $k = 1,2,3,...$, and the set contains frequencies $\Omega_b$ $\omega_k = k\omega_1$, $k = 1,2,3,...$, where β <1 - arbitrary rational number, for example, β = n / m. When d = m and k = n suppression sets will not be empty. It follows that the system of linearly independent components of the observed processes cannot have more than one set containing periodic components that are presented by harmonic series with multiple frequencies. Independence condition components of the measured signals will always be satisfied if the system sets the frequency of random harmonic components will only contain arbitrary incommensurable frequencies.

## 4. The System is Linearly Dependent Set of Signals

It is known[3] that in the Hilbert space of almost periodic functions defined by a completely continuous operator of convolution type of normal form generated by the exact component of the input signal $x(t)$, which shows $B_2(-\infty, +\infty)$ to itself. In multivariate control systems in general, each output additively associated with each entry in a particular control channel. In this case, the exact components of the observed signals synchronously $x_d$ all d inputs and one output r associated operator equation

$$y_r(t) = \sum_{c=1}^{c=d} M\{x_c(t-\tau)k_{rc}(\tau)\} = \sum_{c=1}^{c=d} A_c k_{rc}, \quad (6)$$

where $k_{rc}$- weight function control channel between the output and the input of $r$ $c$.

Instead of the exact component $y_r$ output for each output signal is observed r $\tilde{y}_r$, which belongs to the system of sets of linearly independent components of the output signal $\tilde{y}_r \in \mathbb{Z}_r = Y_r \cup \Gamma_r \cup M_r$, where many elements of accurate output signals $\hat{\theta}_r y_r \in Y_r$, $\gamma_r \in \Gamma_r$_r-set of reduced additive noise generated by the multiplicative interference, $m_r \in M_r$ r-set of additive noise. The integration of these systems, sets all outputs leads to a system of sets of linearly independent components of the output signal $\mathbb{Z} = \cup_{r=1}^{r=h} Z_r$. If the system of sets of input and output signals of $\mathbb{R}$ and $\mathbb{Z}$ are the semicircles, then their union loses this property.

Theorem 2. Suppose that in certain function spaces lot of weight functions $k \in K$ is shown on the set of $Y$ with a linear stationary operator of convolution type $y = Ak$, $y \in Y$, where the kernel of $A$ is generated by the input signal $x \in X$. Image of a compact operator of convolution type is generated[3], spanned by a finite basis orthonormal core, if the operator equation $y = Ak$ defined in the Hilbert space of almost periodic functions in the sense of Besicovitch.

Proof of Theorem 2. It is known[13] that every almost periodic function in the sense of Besicovitch is convergent trigonometric series, for example, $x = \sum_{c=1}^{c=s} a_c b_c$, where $b_c = \exp(j\omega_c t)$. Using the Fourier transform we find that the operator of convolution type is represented in the form $y = \sum_{c=1}^{c=s} \lambda_c \alpha_c b_c$, where $\lambda_c$ - Fourier exponents of the weight function $k \in K$ define the projection of the weight function for the orthonormal function $b_c$. It follows that the range of convolution is a linear span of the orthonormal basis of the kernel $b_c \in G$. It is known[15] that the compact operators generate a finite number of eigenvalues and eigenvectors. Therefore, in the representations of the kernel, image and inverse image Fourier $s < \infty$. Theorem 2.

Consequence of Theorem 2. In the Hilbert space of almost periodic functions is $y \in Y$ and the function that generates the kernel of $x \in X$, a linear stationary operator of convolution type are linearly dependent elements. The linear dependence of the image $y \in Y$ and functions that generate the kernel of $x \in X$, a stationary linear convolution operator follows from the fact that equality $\sum_{c=1}^{c=s} p_c a_c b_c + \sum_{c=1}^{c=s} \lambda_c a_c b_c = 0$ for all nonzero $a_c, b_c, p_c$ and $\lambda_c$, when $p_c = -\lambda_c$. It is known[3] that in the Hilbert space of almost periodic functions $B_2(-\infty, +\infty)$ is an uncountable set of functions $\exp(j\omega t)$, and each vector in this space is only a countable set of non-zero eigenvectors $b_c = \exp(j\omega_c t)$. Theorem 2 implies that the finite set of numeric frequency $\omega_c \in \Omega_x$ and $\omega_c \in \Omega_y$ coincide $\Omega_x = \Omega_y$. It follows that the system of sets ate contains at least one image and inverse image of a linear operator, then the system sets is not a semicircle.

## 5. Ant Blackout with Independent Input Influences

Consider the method of noise control in the special case where all $d$ input control observed linearly independent inputs. In such cases, the system sets the components of input signals $\mathbb{R} = \cup_{k=1}^{k=d} R_k$ is a union of disjoint sets of component signals. By Theorem 1, the system sets the values of the angular frequency of harmonic functions $\mathbb{Q} = \cup \Omega_r$ is a union of disjoint finite sets of numbers. System consolidation sets $\mathbb{Z} = \cup_{r=1}^{r=h} Z_r$, observed on all outputs $h$, is also a semicircle.

Lemma. Let the system of sets of input $\mathbb{R} = \cup_{k=1}^{k=d} R_k$ and output $\mathbb{Z} = \cup_{r=1}^{r=h} Z_r$, the measured signals in multi-dimensional management system containing $d$ inputs and $h$ outputs are semirings. If between the number of input and output number I l there is a linear stationary control channel, represented by the convolution defined in the Hilbert space of almost periodic functions in the sense of Besicovitch, the intersection of number sets angular frequency of orthogonal basis is observed signals at the input and output i l contains finite subset of the numerical values of the angular frequency component of the input signal accurately observed at the input i.

Proof. Semiring of sets $\mathbb{R} = \cup_{k=1}^{k=d} R_k$, semirings is the union of systems of sets, $\mathbb{R}_k$, $k = 1,2,3,..,i,..d$. In this semi-circle is a semicircle $\mathbb{R}_i$ components $\mathbb{R}_i = X_i \cup P_i \cup N_i$, where the infinite set of exact components of the signals $\hat{\theta}x \in X_i$, given additive noise $p \in P_i$, additive noise $n \in N_i$. According to Theorem 1 semicircle of infinite-dimensional sets of elements generated by a semicircle $\mathbb{R}_i$ finite sets of angular frequency harmonic bases

$$\widetilde{\Omega}_{xi} = \Omega_{xi} \cup \Omega_{pi} \cup \Omega_{ni}. \quad (7)$$

By (6) on the output signal is observed l equal to the sum of signals $y_l = \sum_{c=1}^{c=d} A_c k_{cl} = \sum_{c=1}^{c=d} y_{cl}$ where $y_{cl} \in \tilde{Y}_{cl} = \cup_{c=1}^{c=d} Y_{cl} \cup \Gamma_l \cup M_l$, which includes semiring $\tilde{Y}_{il} = Y_{il} \cup \Gamma_l \cup M_l \subset \tilde{Y}_{cl}$. This system is an infinite set of functions generated by a finite number system sets the angular frequency of the harmonic bases $\widetilde{\Omega}_{yil} = \Omega_{yil} \cup \Omega_{\Gamma l} \cup \Omega_{Ml}$, where the finite number sets the angular frequency, respectively, $\Omega_{yil}$ - accurate output signal generated by the operator $A_i k_{il}$ -given additive noise, $\Omega_{\Gamma l}$-additive noise. By Theorem 2 $\Omega_{yil} = \Omega_{xi}$. In this case $\widetilde{\Omega}_{yil} = \Omega_{xi} \cup \Omega_{\Gamma l} \cup \Omega_{Ml}$. We find the suppression of numerical sets of angular frequency basis functions, generating an infinite set of functions on the input i, and number sets the corner frequency of the basis functions that generate an infinite set of functions on the output $l$ $U_{il} = \widetilde{\Omega}_{yil} \cap \widetilde{\Omega}_{xi} = (\Omega_{xi} \cup \Omega_{\Gamma l} \cup \Omega_{Ml}) \cap (\Omega_{xi} \cup \Omega_{pi} \cup \Omega_{ni}) = \Omega_{xi}$. The Lemma is proof.

The lemma is a consequence. Lemma defines the procedure for resolving the problem of isolating the exact components of the signals observed, for example, at the entry and exit of $i$ $l$ in the presence of multiplicative and additive noise. Initially defined numeric arrays and frequency $\widetilde{\Omega}_{yil}$ and $\widetilde{\Omega}_{xi}$, using, for example, the famous ideal analyzer[16] to determine the current amplitude spectrum. Arrays are defined frequencies, for example, the local maximum current amplitude spectra. Then there is the intersection of arrays found frequency $U_{il} = \widetilde{\Omega}_{yil} \cap \widetilde{\Omega}_{xi}$. Allocation projections exact components of the input and output signals observed on the background of multiplicative and additive noise is performed by using the Fourier series of indicator definitions only at frequencies $\omega_k \in U_{il}$.

# 6. Ant Blackout with Dependent Input Actions

In addressing applications as a rule, all or part of the input signals is linearly dependent. In such cases, the property is not satisfied by Lemma and the above described algorithm for applications not applicable. The degree of linear relationship between the input actions can be estimated by the cosine of the angle $\cos(\varphi_{il})$ [14] between the exact components of the signals (prototypes) on the inputs i and $l$. Here are the following cases. If the selected pair of inputs observable transforms signals are linearly independent, then $\cos(\varphi_{il}) = 0$. If the inputs, for example, $i$ and $d$ observed transforms signals coincide up to a sign, then $|\cos(\varphi_{id})| = 1$. The dimension of the control system to the inputs overstated for this case.

One input $i$ or $d$ should be excluded from consideration. After exclusion of extra inputs we find that for systems with linearly dependent input actions for all inputs the condition $|\cos(\varphi_{il})| < 1$.

Theorem 3. Let the system of control has d inputs. Each input signal $\tilde{x}_i$, $1 \leq i \leq d$ is a semi-circle of the signal containing the communication functions (additional prototypes) $f_{il}$, $\tilde{x}_i = \hat{\vartheta}_i x_i + p_i + n_i + \sum_{l=1}^{l=d} f_{il}$ for all $l \neq i$ measure generated by the scalar product in the Hilbert space of almost periodic functions. If all or some of the inputs inputs linear dependence transforms signals $\hat{\vartheta}_i x_i + \sum_{l=1}^{l=d} f_{il}$, satisfying the condition $0 \leq |\cos(\varphi_{ij})| < 1$ for all inputs $i \neq j$, then the original system of linearly dependent input signals $\widetilde{\mathbb{R}} = \cup_{i=1}^{i=d} \widetilde{\mathbb{R}}_i$, containing functions can be identified subsystem linearly independent input signals $\hat{x}_i = \hat{\vartheta}_i x_i + p_i + n_i \in \mathbb{R} \subset \widetilde{\mathbb{R}}$, not containing communications functions.

Proof of Theorem 3. That would qualify for prototypes input $|\cos(\varphi_{il})| < 1$ for all $i \neq l$ represent the measured input signal in the form of $\tilde{x}_i = \hat{\vartheta}_i x_i + p_i + n_i + \sum_{l=1}^{l=d} f_{il}$ for all $i \neq l$. This view differs from the representation (2) the presence of additional prototypes as a liaison between input channels $f_{il}$. By hypothesis, each signal $\tilde{x}_i$ a semicircle components of the measured signals to the selected input. Therefore, all possible scalar products between the constituent signal $\tilde{x}_i$, for example, is the inverse images $(\hat{\vartheta}_i x_i, \bar{f}_{il}) \equiv 0$ identically zero. This condition is satisfied for all $1 \leq i, l \leq d$, $l \neq i$.

Union $\widetilde{\mathbb{R}} = \cup_{i=1}^{i=d} \widetilde{\mathbb{R}}_i$ is a semicircle, as suppression systems of sets input subsystem contains empty sets of communication functions, for example, $\widetilde{\mathbb{R}}_i \cap \widetilde{\mathbb{R}}_l = \sum_{l=1}^{l=d} f_{il}$, $l \neq i$. Let us find easy to implement in practice the symmetric difference[14] between the two systems of sets $\tilde{R}_i \Delta \tilde{R}_l = \widetilde{\mathbb{R}}_i \cup \widetilde{\mathbb{R}}_l \setminus \widetilde{\mathbb{R}}_i \cap \widetilde{\mathbb{R}}_l = \mathbb{R}_i \cup \mathbb{R}_l$. With the available operations, define the union of all symmetric differences for all $i \neq k$ $\cup_{i=1}^{i=d}(\cup_{k=1}^{k=d} \widetilde{\mathbb{R}}_i \Delta \widetilde{\mathbb{R}}_l) = \cup_{i=1}^{i=d} \mathbb{R}_i = \mathbb{R}$ - system sets semirings component inputs, which does not contain a link function. Thus in the original system of linearly dependent input $\widetilde{\mathbb{R}}$, containing a subset of the functions of communication, stands subsystem linearly independent of the input signal $\mathbb{R} \subset \widetilde{\mathbb{R}}$, which does not contain subsets containing communication features $f_{il}$. Theorem 3 is proved.

Theorem 4. Let a linear stationary object control contains $d$ linearly dependent inputs, $h$ linearly dependent outputs and $r \leq dh$ control channels. Each input and each output simultaneously observed images and archetypes, linear time-invariant convolution type operators, distorted linearly independent multiplicative and additive noise defined in the Hilbert space of almost periodic functions of Besicovitch. If among all $d$ prototypes running condition $0 \leq |\cos(\varphi_{ij})| < 1$ for all inputs $i \neq j$, then the frequency method of noise control systems with independent input influence extended to cases where input effects are linearly dependent signals.

Proof of Theorem 4. If the input source control system there is a system of linearly dependent input $\widetilde{\mathbb{R}} = \cup_{i=1}^{i=d} \widetilde{\mathbb{R}}_i$, containing a subset of the functions $f_{il}$ connection between inputs $i$ and $l$, then the randomly selected output $p$ process is observed, generated by the original linearly dependent system of input actions

$$\hat{y}_p(t) = y_{xp}(t) + v_p(t) + \beta_p(t) + \vartheta_p(t), \quad (8)$$

Where $y_{xp}(t) = \sum_{c=1}^{c=d} x_c * k_{cp}$ - the exact component (image) of the output process $x_c(t)$, generated by independent components (prototypes) input signals, a symbol .*. denotes the operator of convolution type, defined in the Hilbert space $B_2(-\infty, +\infty)$.

$v_p(t) = \sum_{r=1}^{r=d-1} \sum_{c=r+1}^{c=d} f_{rc} * (k_{rp} + k_{cp})$ - of the signal generated at the output $p$ of the functions of communication $f_{rc}$ that occurs between the inputs $r$ and $c$, $k_{rp}$, $k_{cp}$, - weighting functions between control channels, respectively, input $r, c$ and output $p$, and $(y_p, \bar{v}_p) = 0$, as both $x_c$ и $f_{rc}$ - are linearly independent random processes for all $c$ and $r$,

$\beta_p(t)$ - this has additive noise generated by the multiplicative noise, which satisfies the linear independence $(\hat{y}_p, \bar{\beta}_p) = (v_p, \bar{\beta}_p) = (\beta_p, \bar{\vartheta}_p) = 0$

$\vartheta_p(t)$ - additive noise at the output of p, satisfying the linear independence

$(\hat{y}_p, \bar{\vartheta}_p) = (v_p, \bar{\vartheta}_p) = (\beta_p, \bar{\vartheta}_p) = 0$.

Implementation processes at different outputs, for example, $p$ and $r$ are linearly dependent processes because of the interconnection between inputs $(\hat{y}_p, \hat{\bar{y}}_r) \neq 0$. However, the components of the process, the observed output $p$ are the sum of mutually orthogonal processes (8). By Theorem 1, the set of frequencies of harmonic components of the process $\hat{y}_p(t)$ is a semicircle of sets of frequencies

$$\widehat{\Omega}_{yp} = \bigcup_{c=1}^{c=d} \Omega_{xc} \bigcup \Omega_{vp} \bigcup \Omega_{\beta p} \bigcup \Omega_{\vartheta p}. \quad (9)$$

Let selection subsystem linearly independent input signals $\tilde{x}_i = \hat{\vartheta} x_{xi} + p_{xi} + n_{xi} \in \mathbb{R}_i \subset \mathbb{R}$, each component of which is a finite Fourier series, for example, (4). Since $\mathbb{R}$ is a subsystem of a semicircle, then by Theorem 1, the system of numerical finite sets of harmonic frequency components of the process of bases $\tilde{x}_i$ is a semicircle with (7). We find the intersection of the selected linearly independent sets of basic harmonic frequency at the input $i$ (7) and the output $p$ (9). Not difficult to see that $\widetilde{\Omega}_{xi} \cap \widehat{\Omega}_{yp} = \Omega_{xi}$ - many basic harmonic frequency components of the signal at the exact $i$ - th input. Knowing a lot of precise frequency components of the input and output signals $\Omega_{xi}$ by using the Fourier transform only the frequencies $\omega_k \in \Omega_{xi}$ find accurate projection components of the input and output signals. Thus, after the separation of independent components of input signals $\mathbb{R} \subset \widetilde{\mathbb{R}}$ of the original system of linearly dependent input signals $\widetilde{\mathbb{R}}$ frequency method of noise control systems with independent input influence extended to cases where input effects are linearly dependent signals. Theorem 4 is proved.

## 7. An Example of Solving the Problem of Identification

The proposed method of noise control is used, for example, in solving the problems of identification of dynamic characteristics of multidimensional control objects in the class of linear stationary models. In particular task was to identify the dynamic characteristics of the aircraft according to the data obtained in single automatic landing. The task of identifying the dynamic characteristics of the above control channels was solved by a new method proposed in the following sequence. First by the local maximum current amplitude spectra[16] of the observed signals were identified frequency source system sets the frequency harmonic input $\widehat{\Omega}_x = \cup_i \widehat{\Omega}_{xi}$, where $\widehat{\Omega}_{xi} = \widetilde{\Omega}_{xi} \cup (\cup_{k=1}^{k=d-1}(\cup_{m=k+1}^{m=d} \Omega_{fkm}))$, $i = 1,2,3,4$.

$\widehat{\Omega}_y$ - set the frequency harmonics output $\hat{y}(t)$, represented by the formula (9).

Given a sign $i$. There are several ways to define the set of frequencies $\widetilde{\Omega}_{xi}$, which does not contain a subset of the frequencies generated by the communication functions $\Omega_{fik}$, for all $k \neq i$. Consider one of them.

Determined by the frequency resolution $\delta = 2\pi/T$, where $T$ -length implementation. The frequencies are considered the same if the absolute difference between two frequencies, taken from the set of frequencies belonging to different inputs, the smaller resolution. For each selected frequency $\omega_{ik} \in \widehat{\Omega}_{xi}$, $k = 1,2,3,...ci$, where $ci$ - the original dimension of the array of frequencies on the selected input $i$. Are the absolute frequency difference $\Delta_{ikmn} = |\omega_{mn} - \omega_{ik}|$ for all frequencies of selected sequences of inputs $m \neq i$. Here, $m,n$ - number of inputs, $n,k$ - number of harmonics, $n = 1,2,3,...cm$, where $cm$ - dimension of the original frequency at the input $m$. If $\Delta_{ikmn} < \delta$, then the system sets the frequency $\widehat{\Omega}_{xi}$ frequency index $ik$ discarded. This gives the system of sets of frequencies $\widetilde{\Omega}_{xi}$, which by (7) does not contain a lot of pairs of matching frequencies generated by communication functions. In this case, the original dimension of the dependent frequency $\widehat{\Omega}_{x2}$ reduced, for example, the original dimension of the array frequency dependent component inlet throttle position was c2= 137, and the dimension of the array of frequencies $\widetilde{\Omega}_{x2}$ independent harmonic o $\tilde{x}_2$ the same second input was c2= 26.

Next, we solve the problem filtering harmonic frequencies useful components of the output signals from the harmonic frequencies of interference. By the corollary of Lemma for the specified control channel $\tilde{x}_2 - y$ defined set of frequencies harmonic generators and describing the forced motion control channel by coincidence frequency harmonic components of the input and output signals of the selected control channel $\widetilde{\Omega}_{x2y} = \{\omega_\rho = \omega_i \in \widetilde{\Omega}_{x2y}: |\omega_i - \omega_k| \leq \delta$ , $\forall \omega i \in \Omega x2, \forall \omega k \in \Omega y$ , $i \in 1, rx2$, $k \in 1, ry$, where $r_y = 35$ - dimension resulting array frequency harmonic output $\Omega_y, \rho \in [1, q]$, where q - the dimension of the same frequency harmonic generators and describing the forced motion in the channel $\tilde{x}_2 - y$, and in this example, it was found that $q = 6$.

For all matching frequencies $\omega_\rho \in \widetilde{\Omega}_{x2y}$ defined input parameters Fourier $S_{x1}(j\omega_\rho) = a(\omega_\rho) + jb(\omega_\rho)$ and output $S_y(j\omega_\rho) = \gamma(\omega_\rho) + j\beta(\omega_\rho)$. At the lowest frequency of multiple frequencies $\omega_{\rho 1} \in \widetilde{\Omega}_{x2y}$ astatism determine the order in the system for the control action. For what is the value of the frequency of the transfer function at the lowest frequency $W(j\omega_{\rho 1}) = S_y(j\omega_{\rho 1})/S_x(j\omega_{\rho 1}) = P + jQ$. In dependence on the position of $W(j\omega_{\rho 1})$ in the complex plane is determined by the order astatism control action. Usually specified point can be located in the first, second or third quadrant of the complex plane, which corresponds to zero, first or second order astatism $p_a$ [8].

These indicators Fourier astatism order (in this example $p_a = 1$) can make a series of algebraic equations from the 2nd to the q - th order, which are the Fourier transforms of linear ordinary differential equations with constant coefficients is unknown from the 2nd to the q - order. For example, the current set of algebraic equations g$\leq q$ -th order is represented as

$\sum_{k=0}^{k=g}(j\omega_{\rho 1})^{k+p_a} T_{k+p_a} S_y(j\omega_{\rho 1}) = S_{x1}(j\omega_{\rho 1})$,
$\sum_{k=0}^{k=g}(j\omega_{\rho c})^{k+p_a} T_{k+p_a} S_y(j\omega_{\rho c}) = S_{x1}(j\omega_{\rho c})$ , where c$\in [2, g]$.

Two systems of equations that this system of equations for the real and imaginary components. Solving these equations sequentially from the 2nd to the q-th order in the coefficients $T_{k+p_a}$, we find that the maximum order of the system of equations $g = q_{max}$, in which the system of equations is not a joint venture. It was found that the dynamic characteristics of the control channel can be described by a linear ordinary differential equation of order 5 with constant coefficients: $T_0 = 0$; $T_1 = -6,6537$; $T_2 = -5,5133$; $T_3 = -3,4828$; $T_4 = -0,6220$; $T_5 = -0,3525$ to a first approximation, for the thrust lever (throttle) - pitch with the elastic deformation of the structure. Coefficients of the differential equation have negative signs under the rules of aerodynamics.[17]. In addition, the differential equation with constant coefficients obtained ninth order for the channel specified - the actual pitch.

## 8. Conclusions

In the theory of stationary random process occurs dualism. The same random stationary process has two incompatible models of representation. Random processes have a continuous spectrum by Wiener in one model. The same processes have a discrete spectrum to the other models in E. Slutsky. It is believed that the first model has a second-order ergodicity, the second such property is not. A widespread model for Wiener thanks ergodic property.

Dualism does not exist in the theory of stationary random processes. Such processes can only have a discrete spectrum. The well-known method of dealing with uncorrelated noise output signals N.Vineru reduced to the use of correlation functions. Correlation functions are not ergodic processes by averaging over the set of realizations. This method of averaging is widely used in theoretical studies, however, found very rarely used in the solution of applied problems. Known statistical dynamics was in a situation where for applications necessary secondary source data in the form of the correlation functions cannot be obtained. In one point of view believe that the random process is an ensemble of realizations, which are the primary source data. Single realization of a random process is a deterministic function. Difficulties in obtaining secondary raw data led to the formulation of the problem of statistical dynamics in the primary source of data in the form of deterministic functions.

The concept of measuring deterministic signal or process is introduced as the sum of three or four linearly independent components. The first component is the exact image or transform operator control channel. The second component is the additive noise. The third component is the reduced additive noise generated by multiplicative noise. This component may be absent in particular cases. The fourth component is a function of the relationship between different input signals, if the multidimensional control system is a control system with dependent inputs. A restriction, all the components of the measured signal describe the motion of material bodies by the second law of Newtonian mechanics.

It is believed that all channels of one-dimensional or multidimensional control systems are described in the first approximation, ordinary differential equations. It follows that the desired function accurate output signal has continuous derivatives up to the n - th order. The exact rules of the differential equation can only be a continuous function. Image and inverse image of the differential operator can be represented by convergent trigonometric series, and the general solution of the differential equation coincides with the forced movement. Such functions cannot be integrated by Riemann and Lebesgue-Stieltjes and integrated over the function, for example, $g(t) = t/2T$. Function with such properties are well-known elements of the Hilbert space of almost periodic functions in the sense of Besicovitch in this space is defined completely continuous normal form of convolution type, which is the inverse of differentiation for the exact components of the measured signals. All components of a measurement signal belong to the system of sets, which is a semi-circle. Measure semiring generated scalar product in the Hilbert space of almost periodic functions $B_2$ the entire line.

In the theory of random processes underlying concept is uncorrelated random processes. In theory nonrandom (deterministic) functions fundamental concept is that semiring system sets of linearly independent components of the signals. It is shown that if an infinite system of the measured signals is a semicircle, the semicircle is a finite set of numerical angular frequency harmonic bases in representations of Fourier components of the measured signals. Frequency harmonic bases image and inverse image of operator equations are the same, so the system sets containing images and preimages is a semicircle. These distinctive properties lead to frequency interference mitigation method, the algorithm is described in the example, where, for example, the control system has four inputs and one output.

At first, we define sets of frequencies of harmonic components of the input signals and output signals. Matching criterion is determined by the frequency resolution. In general, by number input and output number is determined by the selected control channel. The problem of identification is solved sequentially for the selected control channels.

Typically, the system sets the input signals are not a semicircle. It is shown that the system of sets of dependent input signals can be a subsystem of independent sets of input actions. Next, the identification problem is solved by a dedicated subsystem independent input signals. By allocating a subset of matching the frequency of the selected control channel frequency is determined by a subset of the exact harmonic input and output signals. Fourier exponents are found only on a subset of the exact frequency component signals. In this way separates the exact components of the input and output signals from the selected control channel interference. The order and the coefficients of the differential equation found by successive solutions of systems of equations in the frequency domain for increasing orders, while the condition system of equations will be performed. This method was found ordinary differential equation of fifth order for the channel position of levers of the engine - the aircraft pitch and ninth order for the channel - selected - the actual pitch of the aircraft according to the data obtained in single automatic landing with the elastic deformation of the airframe.